# THE ROLE OF INTERPLAY BETWEEN COEFFICIENTS IN THE $G$-CONVERGENCE OF SOME ELLIPTIC EQUATIONS

LUCIO BOCCARDO, LUIGI ORSINA, AND AUGUSTO C. PONCE

*A Ennio De Giorgi, ringraziandolo per il teorema $C^{0,\alpha}$ e per la $\Gamma$-convergenza*

## 1. Introduction and statement of results

Let $\Omega$ be a bounded, open subset of $\mathbb{R}^N$, $\alpha$, $\beta$ in $\mathbb{R}^+$ and $M : \Omega \times \mathbb{R} \to \mathbb{R}^{N^2}$, be a bounded and measurable matrix-valued function such that

$$(1.1) \qquad \alpha |\xi|^2 \leq M(x)\,\xi\,\xi\,, \qquad |M(x)| \leq \beta\,, \quad \text{a.e. } x \in \Omega\,, \ \forall\, \xi \in \mathbb{R}^N\,.$$

In this paper we study the behavior of the solutions $u$ of the linear Dirichlet problems

$$(1.2) \qquad \begin{cases} u \in W_0^{1,2}(\Omega) \cap L^\infty(\Omega) : \\ \displaystyle \int_\Omega M(x)\nabla u \nabla \varphi + \int_\Omega a(x)\,u\varphi = \int_\Omega f(x)\varphi\,, \\ \forall\, \varphi \in W_0^{1,2}(\Omega) \cap L^\infty(\Omega), \end{cases}$$

with respect to perturbations of the matrix $M(x)$ (with respect to the $G$-convergence) and with respect to perturbations of the nonnegative coefficient $a(x)$ and of the right hand side $f(x)$ (with respect to either weak $L^1(\Omega)$ convergence or weak-$*$ convergence as measures) if we assume that

$$(1.3) \qquad \text{there exists } Q > 0 \text{ such that } |f(x)| \leq Q\,a(x) \in L^1(\Omega)\,.$$

The main result in [3] states that, under assumptions (1.1) and (1.3), there exists a weak solution $u$ of (1.2) such that

$$(1.4) \qquad \|u\|_{L^\infty(\Omega)} \leq Q,$$

which also implies

$$\alpha \int_\Omega |\nabla u|^2 \leq Q \int_\Omega a.$$

Thus we have that the interplay between the coefficient of the lower order term and the right hand side yields some regularizing effects on the solution, since $u$ belongs to $W_0^{1,2}(\Omega) \cap L^\infty(\Omega)$ even with a right hand side $f$ belonging only to $L^1(\Omega)$; so that the lower order term is not so "lower" and the principal part is not completely principal.

We recall the definition of $G$-convergence (used in (1.9)) for sequences of elliptic and bounded matrices.

**Definition 1.1.** Let $\{M_n\}$ be a sequence of matrices which satisfies

$$(1.5) \qquad \alpha|\xi|^2 \leq M_n(x)\,\xi\,\xi, \quad |M_n(x)| \leq \beta, \quad \text{a.e. } x \in \Omega, \quad \forall\,\xi \in \mathbb{R}^N\,.$$

The sequence $\{M_n\}$ is said to $G$-converge to a bounded, elliptic matrix $M_0(x)$ if, for every $g$ in $W^{-1,2}(\Omega)$, the sequence $\{w_n\}$ of the unique solutions

$$(1.6) \qquad w_n \in W_0^{1,2}(\Omega), \quad -\mathrm{div}(M_n(x)\nabla w_n) = g \ \text{ in } \Omega,$$

satisfies

$$w_n \rightharpoonup w_0 \text{ weakly in } W_0^{1,2}(\Omega),$$

where $w_0$ is the unique solution of

$$(1.7) \qquad w_0 \in W_0^{1,2}(\Omega), \quad -\mathrm{div}(M_0(x)\nabla w_0) = g \ \text{ in } \Omega.$$

This notion of $G$-convergence was introduced in [17] by S. Spagnolo in the symmetric case. He proved the following compactness theorem: any sequence of symmetric matrices $M_n(x)$ which satisfies (1.5) admits a subsequence which $G$-converges to a matrix $M_0(x)$ of the same type. Then, in [12], E. De Giorgi and S. Spagnolo proved that

$$\lim_{n\to+\infty} \int_\Omega M_n(x)\nabla w_n\nabla w_n = \int_\Omega M_0(x)\nabla w_0\nabla w_0\,.$$

A relationship between $G$-convergence of differential operators and $\Gamma(\text{weak-}W_0^{1,2}(\Omega))$-convergence can be found in [6].

The study in the case of nonsymmetric matrices is due to Murat-Tartar [15] involving the $H$-convergence, where the boundary condition on the Dirichlet problem is removed. Recent contributions can be found in [1] where one finds an alternative proof, using purely variational arguments, of the compactness of $H$-convergence, originally proved by Murat and Tartar. In [2] the authors pursue a variational characterization of the $H$-convergence in terms of the $\Gamma$-convergence of some quadratic forms, including the case of matrices $M_n(x)$ that are possibly nonsymmetric.

In this paper we study the behavior of the sequence $\{u_n\}$ of solutions of the linear Dirichlet problems

$$(1.8) \qquad \begin{cases} u_n \in W_0^{1,2}(\Omega) \cap L^\infty(\Omega) : \\ \displaystyle\int_\Omega M_n(x)\nabla u_n\nabla\varphi + \int_\Omega a_n(x)\,u_n\varphi = \int_\Omega f_n(x)\varphi, \\ \forall\,\varphi \in W_0^{1,2}(\Omega) \cap L^\infty(\Omega), \end{cases}$$

under the assumptions

$$(1.9) \qquad \{M_n\} \ \ G\text{-converges to} \ \ M_0,$$

$$(1.10) \qquad \text{there exists } Q > 0 \text{ such that } |f_n(x)| \leq Q\,a_n(x), \text{where } a_n(x) \in L^1(\Omega), \ \forall\,n;$$

and the sequences

$$(1.11) \qquad \{f_n\} \text{ and } \{a_n\} \text{ converge weakly in } L^1(\Omega) \text{ to } f_0 \text{ and } a_0\,, \text{ respectively.}$$

We emphasize that the constants $\alpha$ and $\beta$ in (1.1) and $Q$ in (1.3) are independent of $n$.

We prove the following result.

**Theorem 1.2.** *Let $\{M_n\}$ be a sequence of matrices satisfying (1.5), which G-converges to a matrix $M_0$. Let $\{a_n\}$ and $\{f_n\}$ be sequences of functions satisfying both (1.10) and (1.11), and let $\{u_n\}$ be the sequence of solutions of (1.8). Then $\{u_n\}$ converges weakly in $W_0^{1,2}(\Omega)$ and weakly-* in $L^\infty(\Omega)$ to the solution $u_0$ of*

$$
\begin{cases}
u_0 \in W_0^{1,2}(\Omega) \cap L^\infty(\Omega) : \\
\displaystyle \int_\Omega M_0(x)\nabla u_0 \nabla \varphi + \int_\Omega a(x)\,u_0\varphi = \int_\Omega f(x)\varphi\,, \\
\quad \forall\,\varphi \in W_0^{1,2}(\Omega) \cap L^\infty(\Omega)\,.
\end{cases}
$$

*Furthermore, we have*

$$
(1.12) \qquad \lim_{n \to +\infty} \int_\Omega M_n(x)\nabla u_n \nabla u_n = \int_\Omega M_0(x)\nabla u_0 \nabla u_0\,.
$$

Since problems (1.8) are stable with respect to both the *G*-convergence of matrices and the weak $L^1$ convergence of the data, one wonders whether this stability is maintained under some weaker convergence of both sequences $\{a_n\}$ and $\{f_n\}$. One could assume for example that both sequences are convergent in the weak-* topology of measures. If this is the case, the main problem one is faced with consists in passing to the limit in integrals of the form

$$
\int_\Omega a_n(x)\,u_n\,\varphi\,,
$$

where $\varphi$ is now a continuous function. An $L^\infty(\Omega)$ estimate on $u_n$ is no longer enough in presence of weak-* convergence of $\{a_n\}$ in the sense of measures: one needs, for example, uniform convergence. A possible approach to prove the uniform convergence of $\{u_n\}$ would be based on some uniform Hölder continuity estimate on $u_n$, which would then allow to prove uniform convergence thanks to the Ascoli-Arzelà theorem. Since for elliptic equations $L^\infty(\Omega)$ results and Hölder continuity results are proved under the same assumptions on the operator and on the data (see [18] and [9]), one may wonder whether also in this case, with data in $L^1(\Omega)$, but solutions in $L^\infty(\Omega)$, some Hölder continuity would hold for solutions of (1.2).

This is not the case: in Section 3 we will give an example in which the functions $a$ and $f$ belong to $L^{N/2}(\Omega)$, and problem (1.2) has a solution which is continuous, and not Hölder continuous. Such a summability is in fact critical: if $a$ and $f$ belong to $L^p(\Omega)$, with $p > \frac{N}{2}$, then any solution $v$ of

$$
-\Delta v + a(x)\,v = f(x)\,,
$$

with $a \geq 0$, is simultaneously bounded by Stampacchia's results (see [18]) and Hölder continuous by De Giorgi's results (see [9]).

The example in Section 3 leaves little hope to prove a stability result in the spirit of Theorem 1.2 using the uniform convergence of the sequence $\{u_n\}$. Actually, the latter property does not hold in general. Nevertheless, we prove in Section 4 that if $\{a_n\}$ and $\{f_n\}$ converge to two bounded Radon measures $\mu$ and $\nu$ such that $|\nu| \leq Q\,\mu$, then the sequence $\{u_n\}$ of solutions of (1.8) converges to the solution $u_0$ of

$$\begin{cases} u_0 \in W_0^{1,2}(\Omega) \cap L^\infty(\Omega): \\ \displaystyle\int_\Omega M_0(x)\nabla u_0 \nabla\varphi + \int_\Omega u_0\,\varphi\,d\mu_d = \int_\Omega \varphi\,d\nu_d\,, \\ \forall\,\varphi \in W_0^{1,2}(\Omega)\cap L^\infty(\Omega)\,, \end{cases}$$

where $\mu_d$ and $\nu_d$ are the absolutely continuous parts of $\mu$ and $\nu$ with respect to the $W^{1,2}$ capacity. This theorem thus yields a result of nonexistence by approximation for (1.2) if $\mu$ and $\nu$ are singular with respect to capacity, as well as a thorough stability result for solutions of (1.8).

## 2. Proof of Theorem 1.2

Before proving Theorem 1.2, we need some technical results; the first one is the Dunford-Pettis theorem.

**Theorem 2.1** (Dunford-Pettis). *Suppose that the sequence $\{y_n\}$ is bounded in $L^1(\Omega)$. Then the sequence is relatively compact in $L^1(\Omega)$ with respect to the weak topology if and only if it is is equi-integrable; that is if for every $\sigma > 0$, there exists $\delta_\sigma > 0$ such that, for any measurable subset $A \subset \Omega$ with $\mathrm{meas}(A) \leq \delta_\sigma$, we have*

$$\sup_{n\in\mathbb{N}} \int_A |y_n| \leq \sigma\,.$$

We will use the following lemma.

**Lemma 2.2.** *Assume that the sequence $\{g_n\}$ converges weakly in $L^1(\Omega)$ to $g_0$ and that the sequence $\{\psi_n\}$ is uniformly bounded and converges almost everywhere to $\psi_0$. Then*

$$(2.1) \qquad\qquad \lim_{n\to+\infty} \int_\Omega g_n\,\psi_n = \int_\Omega g_0\,\psi_0\,.$$

*Proof.* Fix $\sigma > 0$, and let $\delta_\sigma > 0$ be such that, for any measurable subset $A \subset \Omega$ with $\mathrm{meas}(A) \leq \delta_\sigma$, the inequalities

$$(2.2) \qquad\qquad \int_A |g_0| \leq \sigma\,, \qquad \sup_{n\in\mathbb{N}} \int_A |g_n| \leq \sigma$$

hold; use the Dunford-Pettis theorem above and the absolute continuity of the Lebesgue integral to find such a $\delta_\sigma$. We apply Egoroff's theorem in $\Omega$ (which has finite measure, since it is bounded): with $\delta_\sigma > 0$, there exists a measurable subset $\tilde{F}$ such that

$$\mathrm{meas}(\Omega\setminus\tilde{F}) \leq \delta_\sigma\,, \quad \text{and} \quad \psi_n \text{ converges uniformly to } \psi_0 \text{ in } \tilde{F}\,.$$

Taking $M > 0$ such that $\|\psi_n\|_{L^\infty(\Omega)} \leq M$, we have

$$\left| \int_\Omega g_n \psi_n - \int_\Omega g_0 \psi_0 \right| \leq \left| \int_{\tilde F} g_n \psi_n - \int_{\tilde F} g_0 \psi_0 \right| + M \left( \int_{\Omega \setminus \tilde F} |g_n| + \int_{\Omega \setminus \tilde F} |g_0| \right).$$

By uniform convergence of the sequence $\{g_n \psi_n\}$ on $\tilde F$ we have

$$\lim_{n \to +\infty} \int_{\tilde F} g_n \, \psi_n = \int_{\tilde F} g_0 \, \psi_0 \,.$$

From (2.2) with $A = \Omega \setminus \tilde F$, we deduce that

$$\limsup_{n \to +\infty} \left| \int_\Omega g_n \psi_n - \int_\Omega g_0 \psi_0 \right| \leq 2M\sigma \,,$$

which implies (2.1) since $\sigma$ is arbitrary. $\qquad\qquad\qquad\qquad\qquad\qquad\square$

The last tool we need is a consequence of $G$-convergence.

**Lemma 2.3.** *Let $\{M_n\}$ be a sequence of matrices which satisfies (1.5), and that $G$-converges to some matrix $M_0$. Then for every function $\phi_0$ in $W_0^{1,2}(\Omega) \cap L^\infty(\Omega)$, there exists a sequence $\{\phi_n\}$ in $W_0^{1,2}(\Omega) \cap L^\infty(\Omega)$ such that*

$$(2.3) \quad \begin{cases} \{-\mathrm{div}(M_n^*(x)\nabla\phi_n)\} \text{ converges strongly to } -\mathrm{div}(M_0^*(x)\nabla\phi_0) \text{ in } W^{-1,2}(\Omega), \\ \qquad\qquad\qquad |\phi_n| \leq \|\phi_0\|_{L^\infty(\Omega)} \,, \\ \{\phi_n\} \text{ converges to } \phi_0 \text{ weakly in } W_0^{1,2}(\Omega) \text{ and almost everywhere.} \end{cases}$$

*Proof.* Let $z_n$ be the solution of the Dirichlet problem

$$z_n \in W_0^{1,2}(\Omega) : -\mathrm{div}(M_n^*(x)\nabla z_n) = -\mathrm{div}(M_n^*(x)\nabla\phi_0) \,.$$

By Proposition 2 in [15], the sequence $\{M_n^*\}$ of adjoint matrices $G$-converges to $M_0^*$. By definition of $G$-convergence, the sequence $\{z_n\}$ converges weakly in $W_0^{1,2}(\Omega)$ to $\phi_0$. Define, for $k > 0$, the functions $T_k : \mathbb{R} \to \mathbb{R}$ and $G_k : \mathbb{R} \to \mathbb{R}$ by

$$T_k(s) = \begin{cases} -k & \text{if } s \leq -k, \\ s & \text{if } -k < s < k, \\ k & \text{if } s \geq k, \end{cases} \qquad G_k(s) = s - T_k(s) \,.$$

Let $M = \|\phi_0\|_{L^\infty(\Omega)}$, and choose $G_M(z_n)$ as test function in the equation for $z_n$. We have, using (1.5), that

$$\alpha \int_\Omega |\nabla G_M(z_n)|^2 \leq \int_\Omega M_0^*(x)\nabla\phi_0\nabla G_M(z_n) \,.$$

Therefore, since $\{G_M(z_n)\}$ converges weakly in $W_0^{1,2}(\Omega)$ to $G_M(\phi_0) = 0$, we have

$$\limsup_{n \to +\infty} \alpha \int_\Omega |\nabla G_M(z_n)|^2 \leq \lim_{n \to +\infty} \int_\Omega M_0^*(x)\nabla\phi_0\nabla G_M(z_n) = 0 \,,$$

so that $G_M(z_n)$ strongly converges in $W_0^{1,2}(\Omega)$ to $0$. Since $T_M(z_n) = z_n - G_M(z_n)$, the function $\phi_n = T_M(z_n)$ satisfies (2.3). $\qquad\square$

We can now prove Theorem 1.2.

*Proof of Theorem 1.2.* In [3] it is proved that, under the assumption (1.10), there exists a weak solution $u_n$ of (1.8) such that

$$(2.4) \qquad\qquad \|u_n\|_{L^\infty(\Omega)} \le Q \,,$$

and

$$(2.5) \qquad\qquad \alpha \int_\Omega |\nabla u_n|^2 \le Q \int_\Omega a_n \,.$$

The sequence $\{u_n\}$ is bounded in $W_0^{1,2}(\Omega)$; hence, there exists a subsequence (not relabelled) of $\{u_n\}$ that converges weakly in $W_0^{1,2}(\Omega)$ and almost everywhere to some function $u_* \in W_0^{1,2}(\Omega)$. A consequence of (2.4) is that $\|u_*\|_{L^\infty(\Omega)} \le Q$.

Now, fix $\phi_0$ in $W_0^{1,2}(\Omega) \cap L^\infty(\Omega)$, and let $\{\phi_n\}$ be the sequence given by Lemma 2.3 that is contained in $W_0^{1,2}(\Omega) \cap L^\infty(\Omega)$. Using $\phi_n$ as test function in (1.8), we get

$$\int_\Omega M_n^*(x) \nabla \phi_n \nabla u_n + \int_\Omega a_n(x) \, u_n \phi_n = \int_\Omega f_n(x) \phi_n \,.$$

Using the strong convergence in $W^{-1,2}(\Omega)$ of the sequence $\{-\operatorname{div}(M_n^* \nabla \phi_n)\}$ and the weak convergence of $\{u_n\}$ in $W_0^{1,2}(\Omega)$, one has

$$\lim_{n \to +\infty} \int_\Omega M_n^*(x) \nabla \phi_n \nabla u_n = \int_\Omega M_0^*(x) \nabla \phi_0 \nabla u_* \,.$$

On the other hand, Lemma 2.2 applied with $g_n = a_n$ and $\psi_n = u_n \phi_n$, and with $g_n = f_n$ and $\psi_n = \phi_n$, gives

$$(2.6) \qquad \lim_{n \to +\infty} \int_\Omega a_n(x) \, u_n \phi_n = \int_\Omega a_0(x) \, u_* \phi_0 \,, \qquad \lim_{n \to +\infty} \int_\Omega f_n(x) \phi_n = \int_\Omega f_0(x) \phi_0 \,.$$

For later use in the proof we also apply Lemma 2.2, with $g_n = a_n$ and $\psi_n = u_n^2$, and with $g_n = f_n$ and $\psi_n = u_n$, to get

$$(2.7) \qquad \lim_{n \to +\infty} \int_\Omega a_n(x) \, u_n^2 = \int_\Omega a_0(x) \, u_*^2 \,, \qquad \lim_{n \to +\infty} \int_\Omega f_n(x) u_n = \int_\Omega f_0(x) u_* \,.$$

It follows from (2.6) that

$$\int_\Omega M_0(x) \nabla u_* \nabla \phi_0 + \int_\Omega a_0(x) \, u_* \phi_0 = \int_\Omega f_0(x) \phi_0 \,,$$

for every $\phi_0$ in $W_0^{1,2}(\Omega) \cap L^\infty(\Omega)$. Thanks to the uniqueness of solutions for (1.2), proved in [3], we thus have that $u_0 = u_*$. Hence, the whole sequence $\{u_n\}$ converges to $u_0$, as desired.

The uniqueness of the limits in (2.7) also imply that those limits are also true for the entire original sequence $\{u_n\}$. Since we have

$$\int_\Omega M_n(x)\nabla u_n\nabla u_n = \int_\Omega f_n(x)u_n - \int_\Omega a_n(x)\,u_n^2\,,$$

we then deduce that

$$\lim_{n\to+\infty}\int_\Omega M_n(x)\nabla u_n\nabla u_n = \lim_{n\to+\infty}\int_\Omega f_n(x)u_n - \lim_{n\to+\infty}\int_\Omega a_n(x)\,u_n^2$$
$$= \int_\Omega f_0(x)u_0 - \int_\Omega a_0(x)\,u_0^2 = \int_\Omega M_0(x)\nabla u_0\nabla u_0\,,$$

where in the last passage we have used that $u_0$ can be chosen as test function in (1.2) since it belongs to $W_0^{1,2}(\Omega)\cap L^\infty(\Omega)$. This concludes the proof of the theorem. $\qquad\square$

**Remark 2.4.** The convergence result of Theorem 1.2 is related to the $\Gamma(\text{weak-}W_0^{1,2}(\Omega))$-convergence of the functionals

$$J_n(v) = \frac{1}{2}\int_\Omega M_n(x)\nabla v\nabla v + \int_\Omega \Big[\frac{1}{2}a_n(x)v - f_n(x)\Big]v$$
$$= \frac{1}{2}\int_\Omega M_n(x)\nabla v\nabla v + \int_{|v|>2Q}\Big[\frac{1}{2}a_n(x)v - f_n(x)\Big]v + \int_{|v|\le 2Q}\Big[\frac{1}{2}a_n(x)v - f_n(x)\Big]v.$$

Note that the second integral is positive and it can be $+\infty$. Observe that if $M_n$ is symmetric, then the solution $u_n$ of (1.8) is the minimum of $J_n$. The $\Gamma$-convergence (see [10], [11]) is an important tool to prove the convergence of minima of integral functionals. In Theorem 1.2 we proved directly such a convergence.

**Remark 2.5.** In the statement of Theorem 1.2, if instead of (1.9) we make the stronger assumption that the sequence $\{M_n(x)\}$ converges in measure to $M_0(x)$, it is possible to adapt the proof to establish that the sequence $\{u_n\}$ of solutions of (1.8) converges strongly in $W_0^{1,2}(\Omega)$ to the solution $u_0$ of the Dirichlet problem (1.2); such a result is related with the Mosco-convergence (see [14]).

## 3. A counterexample to Hölder continuity

As stated in the Introduction, we now show that, even though solutions of (1.2) are in $L^\infty(\Omega)$, they need not be Hölder continuous in $\Omega$ if $a$ and $f$ are not in $L^p(\Omega)$, with $p > \frac{N}{2}$.

**Proposition 3.1.** *Let $N \ge 2$ and $\Omega = B_1(0)$. There exist nonnegative functions $a(x)$ and $f(x)$ in $L^{\frac{N}{2}}(\Omega)$ satisfying (1.3) and such that problem (1.2) with $M(x) = I$ has a continuous solution which is not Hölder continuous.*

*Proof.* Let $u_\gamma : B_1(0) \to \mathbb{R}$ be defined by $u_\gamma(x) = 1 - |x|^\gamma$ with $0 < \gamma \leq 1$. Since $N \geq 2$, we have $u_\gamma \in W_0^{1,2}(B_1(0))$ and

$$\|\nabla u_\gamma\|_{L^2(B_1(0))} \leq C_1 \,,$$

for some constant $C > 0$ independent of $\gamma$. Then, one has

$$-\Delta u_\gamma(x) = \gamma \, (\gamma - 1)|x|^{\gamma-2} + \gamma \, (N-1)|x|^{\gamma-2} = \gamma \, (N + \gamma - 2)|x|^{\gamma-2} \,.$$

For any $\lambda > 1$, the function

$$a(x) = \frac{1}{|x|^2 \, (\log(2/|x|))^\lambda}$$

belongs to $L^{N/2}(\Omega)$, but not to $L^p(\Omega)$ for every $p > \frac{N}{2}$. Then

$$\frac{-\Delta u_\gamma(x) + a(x) \, u_\gamma(x)}{a(x)} = \gamma \, (N + \gamma - 2)|x|^\gamma (\log{(2/|x|)})^\lambda + 1 - |x|^\gamma =: g_\gamma(x) \,.$$

Observe that

$$0 \leq g_\gamma(x) \leq \frac{C_2}{\gamma^{\lambda-1}} \,,$$

so that $g$ belongs to $L^\infty(\Omega)$, and $\|g\|_{L^\infty(\Omega)} \leq C_2/\gamma^{\lambda-1}$. Thus, if we define

$$f_\gamma(x) = a(x) \, g_\gamma(x) \,,$$

then $f_\gamma \geq 0$ belongs to $L^{N/2}(\Omega)$, and not better, since $g_\gamma(0) = 1$, and $g$ is continuous. Furthermore,

$$(3.1) \qquad\qquad 0 \leq f_\gamma(x) = a(x) \, g_\gamma(x) \leq \frac{C_2}{\gamma^{\lambda-1}} \, a(x) \,.$$

Hence, $u_\gamma$ satisfies the equation

$$-\Delta u_\gamma + a(x) \, u_\gamma = f_\gamma(x) \,,$$

where $a(x)$ and $f_\gamma(x)$ verify property (1.3).

We now choose $\gamma = \gamma_k = \frac{1}{k}$, and define

$$u(x) = \sum_{k=1}^{+\infty} \frac{1}{2^k} \, u_{\gamma_k}(x) \,.$$

This series converges both uniformly and in $W^{1,2}(B_1(0))$. Therefore, $u$ is a continuous function that solves

$$-\Delta u + a(x) \, u = f \,,$$

where

$$f(x) = \sum_{k=1}^{+\infty} \frac{1}{2^k} \, f_{\gamma_k}(x) \,.$$

We now remark that, by estimate (3.1), we have $0 \leq f(x) \leq Q\,a(x)$, where

$$Q = C_2 \sum_{k=1}^{+\infty} \frac{k^{\lambda-1}}{2^k}.$$

On the other hand, the estimate $f \geq f_{\gamma_1}/2$ implies that $f$ is not in $L^p(\Omega)$ for every $p > \frac{N}{2}$. Finally, $u$ is not Hölder continuous; indeed, if $\gamma > 0$ is given and $h$ is a positive integer such that $\frac{1}{h} = \gamma_h < \gamma$, then we have

$$\frac{|u(x) - u(0)|}{|x|^\gamma} = \frac{1 - u(x)}{|x|^\gamma} = \frac{1}{|x|^\gamma} \sum_{k=1}^{+\infty} \frac{1}{2^k}(1 - u_{\gamma_k}(x)) = \frac{1}{|x|^\gamma} \sum_{k=1}^{+\infty} \frac{|x|^{\gamma_k}}{2^k} \geq \frac{|x|^{\gamma_h - \gamma}}{2^h},$$

and the latter quantity diverges as $x$ tends to zero. $\qquad\square$

## 4. Stability of solutions with measure data

The example in the previous section suggests that one should not expect to have uniform convergence of sequences of solutions, which is a useful property to pass to the limit in approximating problems when dealing with measure data. The lack of uniform convergence indeed happens for a precise reason: if one approximates two measures $\mu$ and $\nu$, with $0 \leq \nu \leq Q\,\mu$ for some $Q > 0$, in order to find a solution of

$$-\mathrm{div}(M(x)\,\nabla u) + \mu\,u = \nu\,,$$

then the approximating solutions converge to the solution of *another* problem, and some parts of the measures (the orthogonal parts with respect to $W^{1,2}$ capacity) are lost.

Before stating the precise result, we need some technical tools.

**Lemma 4.1.** *Let $\eta \geq 0$ be a measure in $\mathcal{M}(\Omega)$; then there exist two unique positive measures $\eta_d$ and $\eta_s$ such that*

i) $\eta = \eta_d + \eta_s$;
ii) *$\eta_d$ is absolutely continuous with respect to the $W^{1,2}$ capacity;*
iii) *$\eta_s$ is orthogonal with respect to the $W^{1,2}$ capacity.*

*Furthermore, there exist a function $f$ in $L^1(\Omega)$ and an element $T$ in $W^{-1,2}(\Omega)$ such that $\mu_d = f + T$, in the sense that*

$$\int_\Omega \varphi\,d\eta_d = \int_\Omega f\,\varphi + \langle T, \varphi \rangle\,, \qquad \forall \varphi \in W_0^{1,2}(\Omega) \cap L^\infty(\Omega)\,.$$

*Proof.* We briefly outline the proof of the first part (see e.g. Proposition 14.12 in [16] for the complete argument). The measures $\eta_d$ and $\eta_s$ are obtained by contraction using a Borel set $Z_\eta$ of zero capacity that achieves the supremum

$$\lambda = \sup\{\eta(E)\,,\ \mathrm{cap}(E) = 0\}\,.$$

More precisely, one takes

$$\eta_d = \eta\lfloor_{\Omega \setminus Z_\eta} \quad \text{and} \quad \eta_s = \eta\lfloor_{Z_\eta}.$$

The measure $\eta_s$ is orthogonal with respect to capacity because $\mathrm{cap}(Z_\eta) = 0$, and $\eta_d$ is absolutely continuous by maximality of $Z_\eta$. The uniqueness of the decomposition follows from the observation that if $\eta = \tilde{\eta}_d + \tilde{\eta}_s$ is another decomposition, then the measure $\eta_d - \tilde{\eta}_d = \tilde{\eta}_s - \eta_s$ is simultaneously absolutely continuous and orthogonal with respect to capacity. Hence, it must be identically zero.

The proof of the second part of the result can be found in [5]. $\qquad\square$

**Lemma 4.2.** *Let $0 \le \nu \le Q\,\mu$ be two measures in $\mathcal{M}(\Omega)$. If $\nu = \nu_d + \nu_s$ and $\mu = \mu_d + \mu_s$ are the decompositions of $\nu$ and $\mu$ given by Lemma 4.1, then*
$$0 \le \nu_d \le Q\,\mu_d \quad and \quad 0 \le \nu_s \le Q\,\mu_s\,.$$

*Proof.* Let $Z_\mu \subset \Omega$ be a Borel set of zero capacity such that $\mu_d = \mu\lfloor_{\Omega \setminus Z_\mu}$ and $\mu_s = \mu\lfloor_{Z_\mu}$. By definition, the measure $\nu\lfloor_{Z_\mu}$ is orthogonal with respect to capacity. Since we have
$$0 \le \nu\lfloor_{\Omega \setminus Z_\mu} \le Q\,\mu\lfloor_{\Omega \setminus Z_\mu} = Q\,\mu_d,$$
it follows that the measure $\nu\lfloor_{\Omega \setminus Z_\mu}$ is absolutely continuous with respect to capacity. By the identity $\nu = \nu\lfloor_{\Omega \setminus Z_\mu} + \nu\lfloor_{Z_\mu}$ and the uniqueness of such a decomposition in terms of absolutely continuous and orthogonal parts, we deduce that $\nu_d = \nu\lfloor_{\Omega \setminus Z_\mu}$ and $\nu_s = \nu\lfloor_{Z_\mu}$, and the conclusion follows. $\qquad\square$

We can now state and prove the main result of this section.

**Theorem 4.3.** *Let $\mu \ge 0$ and $\nu \ge 0$ be two measures in $\mathcal{M}(\Omega)$ such that*
$$there\ exists\ Q > 0\ such\ that\ 0 \le \nu \le Q\,\mu.$$
*Let $\{\rho_n\}$ be a sequence of positive $\delta$-approximating convolution kernels, let $\{M_n\}$ be a sequence of matrices which satisfies* (1.5) *and which $G$-converges to a matrix $M_0$, and let $u_n$ in $W_0^{1,2}(\Omega) \cap L^\infty(\Omega)$ be the solution of*
$$\tag{4.1} -\mathrm{div}(M_n(x)\nabla u_n) + (\rho_n * \mu)\,u_n = (\rho_n * \nu)\,.$$
*Then $\{u_n\}$ converges weakly in $W_0^{1,2}(\Omega)$ and weakly-$*$ in $L^\infty(\Omega)$ to the solution $u_0$ in $W_0^{1,2}(\Omega) \cap L^\infty(\Omega)$ of*
$$\tag{4.2} -\mathrm{div}(M_0(x)\nabla u_0) + \mu_d\,u_0 = \nu_d\,,$$
*where $\mu_d$ and $\nu_d$ are the absolutely continuous parts of the measures $\mu$ and $\nu$ with respect to capacity.*

We next explain some tools that are used in the proof of the theorem. The following result is a straightforward consequence of approximation by convolution.

**Lemma 4.4.** *Let $\eta \ge 0$ be a measure in $\mathcal{M}(\Omega)$, decomposed as*
$$\eta = \eta_d + \eta_s = f + T + \eta_s\,,$$
*following the notation of Lemma 4.1. If $\{\rho_n\}$ is a sequence of positive $\delta$-approximating convolution kernels, then*

a) $\rho_n * f \to f$ *strongly in $L^1(\Omega)$;*

b) $\rho_n * T \to T$ strongly $W^{-1,2}(\Omega)$;

c) $\rho_n * \eta_s \to \eta_s$ in the narrow topology of measures.

The approximation by convolution can be nicely paired with suitable convergences (see also Lemma 2.2, where weaker assumptions are made on the sequences involved).

**Lemma 4.5.** *Let $\eta$ be a positive measure in $\mathcal{M}(\Omega)$, decomposed as*

$$\eta = \eta_d + \eta_s = f + T + \eta_s \, .$$

*Let $\{\rho_n\}$ be a sequence of positive $\delta$-approximating convolution kernels, and let $\{\xi_n\}$ be a sequence of functions such that* a') $\xi_n \to \xi$ weakly-$*$ in $L^\infty(\Omega)$;

b') $\xi_n \to \xi$ weakly in $W_0^{1,2}(\Omega)$.

*Then*

$$\lim_{n \to +\infty} \int_\Omega (\rho_n * \eta_d) \, \xi_n = \int_\Omega \xi \, d\eta_d \, .$$

*Proof.* We have

$$\int_\Omega (\rho_n * \eta_d) \, \xi_n = \int_\Omega (\rho_n * f) \, \xi_n + \int_\Omega (\rho_n * T) \, \xi_n \, .$$

The result then follows from items a) and b) of Lemma 4.4, and assumptions a') and b') on $\xi_n$. □

The next result allows to build a family of cut-off functions, starting from sets of zero capacity.

**Lemma 4.6.** *If $\eta \geq 0$ is a measure in $\mathcal{M}(\Omega)$, decomposed as $\eta = \eta_d + \eta_s$, then for every $\sigma > 0$ there exists a function $\psi_\sigma$ in $C_0^\infty(\Omega)$ such that*

i) $0 \leq \psi_\sigma \leq 1$;

ii) $\psi_\sigma \to 0$ weakly-$*$ $L^\infty(\Omega)$ as $\sigma$ tends to zero;

iii) $\psi_\sigma \to 0$ strongly in $W_0^{1,2}(\Omega)$ as $\sigma$ tends to zero;

iv) *it holds*

$$0 \leq \int_\Omega (1 - \psi_\sigma) \, d\eta_s \leq \sigma \, .$$

For the proof of this lemma, we refer the reader to Lemma 5.1 in [8].

*Proof of Theorem 4.3.* Since $0 \leq \rho_n * \nu \leq Q \, \rho_n * \mu$, the solution $u_n$ is bounded in $L^\infty(\Omega)$ by $Q$, which implies that it is bounded in $W_0^{1,2}(\Omega)$. Therefore, up to subsequences, one has

$$u_n \rightharpoonup u_0 \text{ weakly-$*$ in } L^\infty(\Omega), \qquad u_n \rightharpoonup u_0 \text{ weakly in } W_0^{1,2}(\Omega), \qquad u_n \to u_0 \text{ a.e.}$$

Let $\phi_0$ be a function in $W_0^{1,2}(\Omega) \cap L^\infty(\Omega)$, and let $\{\phi_n\}$ be the sequence of functions in $W_0^{1,2}(\Omega) \cap L^\infty(\Omega)$ given by Lemma 2.3. Let $\psi_\sigma$ be the function given by Lemma 4.6 for $\mu_s$ and $\sigma > 0$. Choosing $\phi_n(1 - \psi_\sigma)$ as test function in (4.1), we have

(4.3)
$$\int_\Omega M_n^*(x)\nabla\phi_n\nabla\big[u_n(1-\psi_\sigma)\big] + \int_\Omega (\rho_n * \mu_d)\, u_n\,\phi_n\,(1-\psi_\sigma) - \int_\Omega (\rho_n * \nu_d)\,\phi_n\,(1-\psi_\sigma)$$
$$= -\int_\Omega (\rho_n * \mu_s)\, u_n\,\phi_n\,(1-\psi_\sigma) + \int_\Omega (\rho_n * \nu_s)\,\phi_n\,(1-\psi_\sigma).$$

We now have, by Lemma 2.3, and by Lemma 4.5 applied once with $\eta = \mu$ and $\xi_n = u_n\,\phi_n\,(1-\psi_\sigma)$, and once with $\eta = \nu$ and $\xi_n = \phi_n\,(1-\psi_\sigma)$, that

$$\int_\Omega M_n^*(x)\nabla\phi_n\nabla\big[u_n(1-\psi_\sigma)\big] + \int_\Omega (\rho_n * \mu_d)\, u_n\,\phi_n\,(1-\psi_\sigma) - \int_\Omega (\rho_n * \nu_d)\,\phi_n\,(1-\psi_\sigma)$$
$$\longrightarrow \int_\Omega M_0^*(x)\nabla\phi_0\nabla\big[u_0\,(1-\psi_\sigma)\big] + \int_\Omega u_0\,\phi_0\,(1-\psi_\sigma)\,d\mu_d - \int_\Omega \phi_0\,(1-\psi_\sigma)\,d\nu_d\,,$$

as $n \to +\infty$. On the other hand, since $M_0^*\nabla\phi_0$ belongs to $(L^\infty(\Omega))^N$, and $u_0\,(1-\psi_\sigma)$ converges strongly in $W_0^{1,2}(\Omega)$ to $u_0$, we have

$$\lim_{\sigma \to 0} \int_\Omega M_0^*(x)\nabla\phi_0\nabla\big[u_0\,(1-\psi_\sigma)\big] = \int_\Omega M_0^*(x)\nabla\phi_0\nabla u_0 = \int_\Omega M_0(x)\nabla u_0\nabla\phi_0.$$

Recall that $\psi_\sigma$ converges to 0 both weakly-$*$ in $L^\infty(\Omega)$ and strongly in $W^{1,2}(\Omega)$. Thus the same holds for $u_0\,\phi_0\,\psi_\sigma$ and $\phi_0\,\psi_\sigma$. Since by Lemma 4.1 the measures $\mu_d$ and $\nu_d$ can be written as a sum of elements in $L^1(\Omega)$ and in $W^{-1,2}(\Omega)$, we have

$$\lim_{\sigma \to 0} \int_\Omega u_0\,\phi_0\,\psi_\sigma\,d\mu_d = 0 \quad \text{and} \quad \lim_{\sigma \to 0} \int_\Omega \phi_0\,\psi_\sigma\,d\nu_d = 0.$$

Therefore,

$$\int_\Omega M_0^*(x)\nabla\phi_0\nabla\big[u_0(1-\psi_\sigma)\big] + \int_\Omega u_0\,\phi_0\,(1-\psi_\sigma)\,d\mu_d - \int_\Omega \phi_0\,(1-\psi_\sigma)\,d\nu_d$$
$$\longrightarrow \int_\Omega M_0(x)\nabla u_0\nabla\phi_0 + \int_\Omega u_0\,\phi_0\,d\mu_d - \int_\Omega \phi_0\,d\nu_d\,,$$

as $\sigma \to 0$.

Observe now that by Lemma 4.2 we have $0 \le \nu_s \le Q\,\mu_s$, so that $\rho_n * \nu_s \le Q\,\rho_n * \mu_s$. Hence,

$$\left| \int_\Omega (\rho_n * \nu_s)\,\phi_n\,(1-\psi_\sigma) \right| \le \|\phi_n\|_{L^\infty(\Omega)} \int_\Omega (\rho_n * \nu_s)\,(1-\psi_\sigma) \le Q\,\|\phi_0\|_{L^\infty(\Omega)} \int_\Omega (\rho_n * \mu_s)\,(1-\psi_\sigma).$$

Similarly, since $\|u_n\|_{L^\infty(\Omega)} \le Q$, we have

$$\left| \int_\Omega (\rho_n * \nu_s)\, u_n\,\phi_n\,(1-\psi_\sigma) \right| \le Q\,\|\phi_0\|_{L^\infty(\Omega)} \int_\Omega (\rho_n * \mu_s)\,(1-\psi_\sigma).$$

Thus,

$$\left| \int_\Omega (\rho_n * \mu_s) \, u_n \, \phi_n \, (1 - \psi_\sigma) \right| + \left| \int_\Omega (\rho_n * \nu_s) \, \phi_n \, (1 - \psi_\sigma) \right| \leq 2Q \, \|\phi_0\|_{L^\infty(\Omega)} \, \sigma \,,$$

which implies that

$$\lim_{\sigma \to 0} \lim_{n \to \infty} \left| \int_\Omega (\rho_n * \mu_s) \, u_n \, \phi_n \, (1 - \psi_\sigma) \right| + \left| \int_\Omega (\rho_n * \nu_s) \, \phi_n \, (1 - \psi_\sigma) \right| = 0 \,.$$

Thus letting first $n \to +\infty$ and then $\sigma \to 0$ in (4.3), we deduce that $u_0$ in $W_0^{1,2}(\Omega) \cap L^\infty(\Omega)$ satisfies

$$\int_\Omega M_0(x) \nabla u_0 \nabla \phi_0 + \int_\Omega u_0 \, \phi_0 \, d\mu_d - \int_\Omega \phi_0 \, d\nu_d = 0 \,,$$

for every $\phi_0$ in $W_0^{1,2}(\Omega) \cap L^\infty(\Omega)$; i.e., $u_0$ is the solution of (4.2). By uniqueness of the solution, the whole sequence $\{u_n\}$ converges to $u_0$. □

**Remark 4.7.** The previous result states that if $\mu = \mu_s$ is orthogonal to capacity (so that $\nu = \nu_s$ is orthogonal to capacity as well), then the sequence $\{u_n\}$ of solutions of (4.1) tends to zero; i.e., there is no solution obtained by approximation for the limit problem

$$-\mathrm{div}(M_0(x) \nabla u) + \mu \, u = \nu \,.$$

This is mainly due to the assumption $0 \leq \nu \leq Q \, \mu$, which yields bounded solutions in $W_0^{1,2}(\Omega)$. Indeed, if such an assumption is missing (so that $\nu$ is not related to $\mu$), then a solution of

$$-\mathrm{div}(M_0(x) \nabla u) + \mu \, u = \nu \,,$$

always exists, provided $\mu$ is absolutely continuous with respect to capacity, and $\nu$ is any bounded measure; see [13] and [7]. In general, the solution in this case does not belong to $W_0^{1,2}(\Omega)$, and can be found by duality techniques in the larger space $W_0^{1,q}(\Omega)$, for every $q < \frac{N}{N-1}$.

**Remark 4.8.** As a consequence of Theorem 4.3, it is possible to give a negative answer to a question raised by Piero Marcati to the first author in a personal communication: if $\{h_n\}$ is a sequence of positive functions converging to $\delta_0$, the Dirac mass concentrated at the origin, and if $\{u_n\}$ is the sequence of solutions of

$$u_n \in W_0^{1,2}(\Omega) \cap L^\infty(\Omega): \ -\mathrm{div}(M(x) \nabla u_n) + h_n \, u_n = h_n \,,$$

then $\{u_n\}$ converges weakly in $W_0^{1,2}(\Omega)$ to zero. Indeed, in this case $\mu = \nu = \delta_0$ are orthogonal with respect to capacity in dimension $N \geq 2$. Thus, the sequence $\{u_n\}$ converges to the unique solution $u_0$ of $-\mathrm{div}(M(x) \nabla u_0) = 0$, which is zero.

This fact is not surprising; indeed, if $u$ is a (continuous) solution of

$$-\mathrm{div}(M(x) \nabla u) + u(0) \, \delta_0 = \delta_0 \,,$$

then $0 \leq u \leq 1$. If $u(0) < 1$, then the function $u$ is a solution of $-\mathrm{div}(M(x) \nabla u) = (1 - u(0)) \delta_0$. Solutions of this equation are unbounded at the origin in dimension $N \geq 2$, which yields a contradiction with the estimate $0 \leq u \leq 1$. On the other hand, if $u(0) = 1$, then $-\mathrm{div}(M(x)\nabla u) = 0$, hence $u \equiv 0$, and we again reach a contradiction.

If, instead, $\mu$ and $\nu$ are functions in Lebesgue spaces, the situation is rather different. Indeed, if $a(x) = \frac{1}{|x|^q}$, with $2 \leq q < N$, then $u(x) = 1 - |x|^\gamma$, with $\gamma > 0$, is a solution of

$$-\Delta u + a(x)\, u = \frac{1}{|x|^q} - \frac{1}{|x|^{q-\gamma}} + \gamma(N + \gamma - 2)\, \frac{1}{|x|^{2-\gamma}} = f(x)\,.$$

Both $f$ and $a$ do not belong to $L^{N/2}(\Omega)$, and that $|f(x)| \leq C_\gamma a(x)$ for some constant $C_\gamma$, so that assumption (1.3) is satisfied. In this case, $u(0) = 1$ does not yield any contradiction. Observe that in this case the function

$$f(x) - a(x)\, u = \gamma(N + \gamma - 2)\, \frac{1}{|x|^{2-\gamma}}\,,$$

belongs to $L^p(\Omega)$ for some $p > \frac{N}{2}$.

**Remark 4.9.** Since the equation in (4.2) is linear, the conclusion of Theorem 4.3 also holds if $\nu$ is a signed measure such that $|\nu| \leq Q\,\mu$ for some $Q > 0$.

**Remark 4.10.** Theorem 4.3 also has a counterpart for sequences of functions $\{a_n\}$ and $\{f_n\}$ which are such that $|f_n| \leq Q\,a_n$, and satisfy the following assumptions: if $\mu = g + T + \mu_s$ is decomposed as in Lemma 4.1, then $a_n = a_{n,1} + a_{n,2} + a_{n,3}$, with $a_{n,1}$ converging to $g$ weakly in $L^1(\Omega)$, $a_{n,2}$ converging to $T$ strongly in $W^{-1,2}(\Omega)$, and $a_{n,3} \geq 0$ converging to $\mu_s$ in the narrow topology of measures. Furthermore, if $\nu = f + S + \nu_s$ is decomposed as in Lemma 4.1, then $f_n = f_{n,1} + f_{n,2} + f_{n,3}$, with $f_{n,1}$ converging to $f$ weakly in $L^1(\Omega)$, $f_{n,2}$ converging to $S$ strongly in $W^{-1,2}(\Omega)$, and $f_{n,3} \geq 0$ converging to $\nu_s$ in the narrow topology of measures.

## Acknowledgements

The third author (ACP) was supported by the Fonds de la Recherche scientifique – FNRS under research grant J.0026.15. He warmly thanks the Dipartimento di Matematica of the "Sapienza" Università di Roma for the invitation. He also acknowledges the hospitality of the Academia Belgica in Rome.

Lucio Boccardo
"Sapienza" Università di Roma
Dipartimento di Matematica
P.le A. Moro 2
00185 Roma
Italy

Luigi Orsina
"Sapienza" Università di Roma
Dipartimento di Matematica
P.le A. Moro 2
00185 Roma
Italy

Augusto C. Ponce
Université catholique de Louvain
Institut de Recherche en Mathématique et Physique
Chemin du cyclotron 2, bte L7.01.02
1348 Louvain-la-Neuve
Belgium